\documentclass[12pt]{amsart}

\newtheorem{theorem}{Theorem}[section]
\newtheorem{lemma}[theorem]{Lemma}
\newtheorem{claim}[theorem]{Claim}
\newtheorem{corollary}[theorem]{Corollary}
\newtheorem{proposition}[theorem]{Proposition}

\theoremstyle{definition}
\newtheorem{definition}[theorem]{Definition}
\newtheorem{question}[theorem]{Question}
\newtheorem{example}[theorem]{Example}

\theoremstyle{remark}
\newtheorem{remark}[theorem]{Remark}

\numberwithin{equation}{section}

\newcommand{\Hh}{\mathbb{H}}
\newcommand{\D}{\mathbb{D}}
\newcommand{\C}{\mathbb{C}}

\newcommand{\St}{\mathbb{S}}
\newcommand{\R}{\mathbb{R}}

\newcommand{\ep}{\epsilon }

\newcommand{\de}{\delta }
\newcommand{\De}{\Delta }

\newcommand{\ga}{\gamma }

\newcommand{\Om}{\Omega }

\newcommand{\rea}{\operatorname{Re}}
\newcommand{\ima}{\operatorname{Im}}
\newcommand{\Arg}{\operatorname{Arg}}

\newcommand{\latlim}{\operatorname{lat.-lim}}

\newcommand{\bd}[1]{\partial #1}

\begin{document}
\baselineskip=18pt

\title[Backward iteration]{Backward-iteration sequences with bounded
hyperbolic steps for analytic self-maps of the disk}  
\author{Pietro Poggi-Corradini}

\date{July 26, 2001}
\address{Department of Mathematics, Cardwell Hall, Kansas State University,
Manhattan, KS 66506}
\email{pietro@math.ksu.edu}



\begin{abstract}
A lot is known about the forward iterates of an analytic function
which is bounded by $1$ in modulus on
the unit disk $\D$. The
Denjoy-Wolff Theorem describes their convergence properties and several
authors, from the 1880's to the 1980's, have provided conjugations which
yield very precise descriptions of the dynamics. Backward-iteration
sequences are of a different nature because a point could have
infinitely many preimages as well as none. 
However, if we insist in choosing preimages
that are at a finite hyperbolic distance each time, we obtain
sequences which have many similarities with the forward-iteration
sequences, and which also reveal more information about the map itself.
In this note we try to present a complete study of backward-iteration
sequences with bounded hyperbolic steps for analytic self-maps of the disk.
\end{abstract}

\maketitle

\section{Introduction}\label{sec:intro}
Let $\phi $ be an analytic self-map
of $\D$. The Theorem of Denjoy-Wolff says that, aside for the case
when $\phi $ is an elliptic automorphism, there is a point 
$\tau_{\phi}\in \overline {\D}$ (which we call the Denjoy-Wolff
point of $\phi$) such that the iterates of $\phi $ converge to
$\tau_{\phi }$ uniformly on compact subsets of $\D$. A lot is known
about the behavior of a forward-iteration sequence under $\phi$,
$z_{n}=\phi _{n} (z_{0})$, especially in the limit, for large values of
$n$. In fact, conjugations of $\phi $ to linear maps have been
established near the Denjoy-Wolff point and these provide a fairly
clear picture of the dynamics. The different types of  behavior divide the
class of self-maps of the disk into three main cathegories:\\
\noindent  {\bf Elliptic.} The Denjoy-Wolff point 
is an interior 
fixed point for $\phi$, i.e. $\tau_{\phi }\in \D$. By Schwarz's Lemma,
the derivative $\lambda $ of $\phi$ at 
$\tau_{\phi }$ (also known as the {\sf muliplier}) satisfies
$|\lambda|\leq 1$. If $|\lambda |=1$, $\phi $ is an
elliptic automorphism. If $0<|\lambda |<1$, $\phi (z) $ can
be conjugated to $\lambda z$ near $\tau _{\phi }$, 
in such a way that the forward
orbit $z_{n}$ is asymptotic to $\lambda^{n}\omega$ for some $\omega$, for large $n$. 
If $\lambda =0$,  $\phi (z) $ can
be conjugated to $z^{N}$ near $\tau_{\phi }$, for some $N\geq 2$, and
related asymptotics can be obtained.\\
\noindent {\bf Hyperbolic.} The Denjoy-Wolff point is a boundary fix
point for $\phi$, i.e., $\tau_{\phi }\in \bd 
\D$ and $\phi (\tau_{\phi
})=\tau_{\phi }$ in the sense of non-tangential limits. In this case, the
derivative $c$ of $\phi$ at $\tau_{\phi }$ exists, again  in the sense 
of non-tangential limits, and satisfies
$0<c\leq 1$. The map 
$\phi$ is hyperbolic if $c<1$.
Conjugations exist in this case and can be used to show that
$\Arg (z_{n}-\tau _{\phi })-\Arg \tau _{\phi } $ tends to an angle
$\theta _{0}\in (-\pi /2,\pi /2)$. 
In particular, a forward-iteration sequence eventually tends to
$\tau_{\phi }$ along a non-tangential ray at $\tau_{\phi }$ (Combine
Lemma 2.66 of \cite{cm} with Theorem 3 (i) of \cite{pom}). \\
\noindent  {\bf Parabolic.} The Denjoy-Wolff point $\tau_{\phi }\in \bd
\D$ is a boundary fixed point for $\phi$, and the
derivative $c$ of $\phi$ at $\tau_{\phi }$ is equal to $1$.
This case is more subtle and the behavior of forward-iteration
sequences has been studied by Pommerenke and Baker-Pommerenke in
\cite{pom} and \cite{bp}, and by Cowen, see \cite{cm}. 
Conjugations still exist:
they send $\tau_{\phi }$ to $\infty $ and conjugate $\phi (z)$ to
a translation, see Section \ref{ssec:forward} for more details.

Let $d$ be the pseudo-hyperbolic distance in $\D$, i.e., $d
(z,w)=|z-w|/|1-\overline{z}w|$.   
By Schwarz-Pick, forward-iteration sequences always have bounded steps
in the hyperbolic metric. More specifically, $d (z_{n+1},z_{n})\leq d
(z_{n},z_{n-1})$, and this is maybe the reason why these sequences
become more and more ``regular'' in the limit. Actually, the sequence
of step-lengths, $s_{n}=d (z_{n+1},z_{n})$, plays an important
role. Clearly $s_{n}$ tends to $s_{\infty }\geq 0$. In the elliptic
case $s_{\infty }=0$ (aside for elliptic automorphisms). In the
hyperbolic case, $s_{\infty }>0$. More importantly, in the parabolic
case, both $s_{\infty }=0$ and $s_{\infty }>0$ can occur, and the
dynamics is different in each case. As a matter of notation, we say
$\phi $ is {\bf type I parabolic} if $s_{\infty }>0$, and {\bf type II
parabolic} if $s_{\infty }=0$.

\subsection{Boundary repelling fixed points}\label{ssec:brfps}

We saw above that $\tau_{\phi }$ is
always a fixed point for $\phi $ and the multiplier there is less or
equal to $1$ in modulus. The Denjoy-Wolff Theorem says, moreover, that
every other fixed point of $\phi $ can only be on the boundary of $\D$,
with positive multiplier strictly greater than $1$ (infinity an
admitted value). An analytic self-map of
the disk may have several such fixed points other than the Denjoy-Wolff
point, however, when the
multiplier is also finite then more can be said and these points
have several nice features. That is the reason why in \cite{pams} we
gave them the name of Boundary Repelling Fixed Points (BRFP).
Below are three results involving boundary fixed points which will be useful in
the sequel. 

\begin{theorem}[Julia and Carath\'eodory,\cite{sh} Chap. 4]\label{thm:juca}
Suppose $\phi$ is an analytic map of the disk with
$\phi(\D)\subset\D$, and $\zeta,\xi \in\bd \D$. 
If there is a sequence
$\{p_{n} \}\subset \D$ such that $p_{n}\rightarrow \zeta$, $\phi
(p_{n})\rightarrow \xi$, and 
\begin{equation}\label{eq:julia}
\frac{1-|\phi(p_{n})|}{1-|p_{n}|}\rightarrow A <\infty,
\end{equation}
then
\begin{enumerate}
\item [(a)] $A>0$
\item [(b)] For every horodisk $H$ at $\zeta $, i.e., $H$ is a disk
internally tangent to $\bd\D$ at $\zeta $, $\phi (H)\subset M
(H)$, where $M (z)=\xi\overline{\zeta }(z-a\zeta)/(1-a\overline
{\zeta }z)$, with $a=\frac{A -1}{A+1}$.
\item [(c)] $\phi (z)\rightarrow \xi$ as $z\rightarrow \zeta $
non-tangentially.     
\item [(d)] $\phi ^{\prime } (z)\rightarrow \phi ^{\prime } (\zeta )$
as $z\rightarrow \zeta $ non-tangentially, and $|\phi^{\prime } (\zeta
)|\leq A$.
\end{enumerate}
\end{theorem}
\begin{corollary}\label{thm:jc}
Suppose $\phi$ is an analytic map with
$\phi(\D)\subset\D$, and suppose $\zeta \in\bd \D$. Assume further that
there is a sequence
$\{p_{n} \}\subset \D$ such that 
\begin{enumerate}
\item $p_{n}\rightarrow \zeta$,
\item $\lim_{n\rightarrow \infty } d(p_{n},\phi (p_{n}))\leq a<1$.
\end{enumerate}
Then $\zeta $ is a boundary fixed point of $\phi$ with multiplier $\phi ^{\prime } (\zeta )\leq
\frac{1+a}{1-a}$. 
\end{corollary}

\begin{theorem}[Cowen-Pommerenke, see Theorem 4.1 of \cite{cp}]\label{thm:cp}
Let $\phi$ be analytic with $\phi(\D)\subset\D$, and let $\tau_{\phi}$ be its
Denjoy-Wolff point. If $\tau _{\phi }\in \D$ or if $\tau _{\phi }\in
\bd \D$ and $\phi ^{\prime } (\tau _{\phi })<1$ (elliptic and
hyperbolic cases), then for every $A>1$ the set of BRFPs whose
multiplier is less than $A$ must be finite. Moreover if $\tau _{\phi }\in
\bd \D$ and $\phi ^{\prime } (\tau _{\phi })=1$ (parabolic case), then
for every $A>1$ the set of BRFPs whose multiplier is less than $A$ can
only cluster at $\tau _{\phi }$.
\end{theorem}
The two theorems just cited were used in \cite{finn} to prove the
following conjugation result.
\begin{theorem}[\cite{finn} Thm. 1.2]\label{thm:finn}
Suppose $\phi$ is analytic with $\phi(\D)\subset\D$.
Assume that $1$ is a BRFP for $\phi$ with multiplier $1<A<\infty$.

Then there is an analytic map  
$\psi$ of the upper half-plane $\Hh$, with $\psi (\Hh)\subset \D$, 
which has non-tangential limit $1$ at $0$, and such that:
\begin{equation}\label{eq:conj}
\psi(Az)=\phi\circ\psi(z)
\end{equation}
for every $z\in\Hh$. 

Moreover, $\psi$ is always semi-conformal at $0$, 
that is, 
\begin{equation}\label{eq:semiconf}
\Arg\frac{\psi(z)-1}{iz}\longrightarrow 0
\end{equation}
as $z$ tends to $0$ non-tangentially.
\end{theorem}

The main construction in the proof of Theorem \ref{thm:finn} was to
produce a backward-iteration sequence with certain additional
properties, in particular whose steps remain bounded in the
hyperbolic distance. 
\begin{definition}
A sequence $\{w_{n}\}_{0}^{\infty }$  is a
backward-iteration sequence with bounded steps (BISBS) for $\phi$, if $\phi
(w_{n+1})=w_{n}$ for $n=0,1,2,3\dots $, and $d (w_{n},w_{n+1})\leq
a<1$ for all $n$ and for some constant $a$. 
(We exclude the trivial
sequence $w_{n}\equiv \tau_{\phi }$, in the 
elliptic case, from this definition, and we also assume that $\phi $ is
not an elliptic automorphism in the rest of this paper.)
\end{definition}
On the other hand, Theorem \ref{thm:finn} implies
that if $\phi $ has a BRFP at $\zeta $ (say $\zeta =1$) then there are a lot of
BISBS, namely all the sequences of the
form $\{\psi (A^{-n}z) \}_{n=0}^{\infty }$ for some $z\in\Hh $.

In this note, we propose to study the BISBS of a self-map $\phi $ of
the disk, their convergence
properties, as well as existence and uniqueness properties.
To begin with, in the result stated in the next section we show that a
BISBS can only arise as in Theorem \ref{thm:finn} above, i.e. tends to a BRFP non-tangentially and is 
of the form $\{\psi (A^{-n}z) \}_{n=0}^{\infty }$ for a conjugation as
in Theorem \ref{thm:finn}, except for the parabolic case, where it can
happen that a BISBS
actually tends to the Denjoy-Wolff point. The parabolic case is more
delicate and is the main focus of this paper, see Paragraphs
\ref{ssec:forward},\ref{ssec:backward}, \ref{ssec:uquests},
\ref{ssec:conjtyp1}, \ref{ssec:nonexist}   below.

\subsection{Backward-iteration sequences with bounded
steps}\label{ssec:bisbs}

 We will
see that the bounded-steps
restriction yields several interesting properties and  that
these sequences become ``regular'' for large values of $n$, 
analoguously to the forward-iteration sequences.

Note that Schwarz-Pick implies that $d (w_{n},w_{n+1})$ is increasing
with $n$, thus without loss of generality we always assume that
$d (w_{n},w_{n+1})\uparrow a<1$.
\begin{theorem}\label{thm:converges}
Suppose $\phi $ is an analytic map with $\phi(\D)\subset\D$, and let
$\{w_{n}\}_{n=0}^{\infty }$ be a backward-iteration sequence for $\phi
$ with bounded pseudo-hyperbolic steps $d_{n}=d (w_{n},w_{n+1})\uparrow 
a<1$. Then, the following
hold:
\begin{enumerate}
\item There is a point
 $\zeta \in \bd  \D$ such that $w_{n}\rightarrow \zeta $ as
$n$ tends to infinity, and $\zeta $ is a fixed point for $\phi $
with a well-defined  multiplier $\phi ^{\prime } (\zeta )<\infty$.
\item When $\zeta\neq \tau _{\phi } $, then
$\zeta $ is a
BRFP. If $\zeta=\tau_{\phi}$, then $\phi $ is necessarily of parabolic
type.
\item When $\zeta $ is a BRFP, 
the sequence $w_{n}$ tends to $\zeta $ along a 
non-tangential direction. More precisely, there exists an angle
$\theta_{0}$ in $(-\pi
/2,\pi /2)$ such that 
\begin{equation}\label{eq:angappr}
\arg (\zeta -w_{n})- \arg \zeta \rightarrow \theta _{0}
\end{equation}
as $n\rightarrow \infty $.
\item When, in the parabolic case, $\zeta =\tau_{\phi }$, then  $w_{n}$
tends to $\zeta $ tangentially.   
\end{enumerate}
\end{theorem}
As we have seen above, 
a consequence of Theorem \ref{thm:finn} is that whenever $\zeta
\in \bd  \D$ is a BRFP for $\phi $ then one can construct many BISBS
converging to $\zeta $ along non-tangential directions as in
(\ref{eq:angappr}). One of the consequences of Theorem
\ref{thm:converges} (3) is that every BISBS approaching a BRFP must do
so non-tangentially as in (\ref{eq:angappr}). 
Finally, as mentioned in (2), in the parabolic case it
can happen that a BISBS  tends to $\tau_{\phi } (=1)$, and we can
construct examples in both the type I and type II cases, see Section
\ref{sec:examples}. Theorem \ref{thm:converges} is proved in Section
\ref{sec:convergence}. The rest of the paper is devoted to studying what
happens in the parabolic case.

\subsection{Parabolic forward iteration}\label{ssec:forward}

We study the parabolic case in Section \ref{sec:parabolic}, where we
obtain counterparts to results about forward iterates of Pommerenke
and Baker-Pommerenke, \cite{pom} \cite{bp}.
First let us recall what is known about forward-iteration in this case.
A self-map of $\D$ is parabolic if $\tau_{\phi }\in\bd \D$ and $\phi
^{\prime } (\tau _{\phi })=1$ there.
Without loss of generality, we assume instead that $\phi $ is an analytic map
on the upper half-plane $\Hh $, with $\ima \phi (z)\geq 0$, such that 
\begin{equation}\label{eq:assump}
\phi (z)\longrightarrow \infty \qquad \mbox { and
}\qquad \frac{\phi (z)}{z}\longrightarrow 1.
\end{equation}
as $z\rightarrow \infty $ non-tangentially.
The non-tangential approach regions for $\infty $ in $\Hh $ are the
sectors $\{|\Arg z-\pi /2|<\theta_{0}  \}$ for some $\theta_{0} \in
(0,\pi /2)$. Since the horodisks at infinity are the half-planes
$\{\ima z\geq t>0 \}$, Julia's Lemma in this situation implies that
$\ima \phi (z)\geq \ima z$. 
In $\Hh $ the pseudo-hyperbolic distance between two
points $z,w$ is:
\[
d=d (z,w)=\left|\frac{w-z}{w-\overline {z}} \right|
\]
Let $z_{n}=\phi _{n} (i)$. We know that $z_{n}$ tends to infinity by
the Denjoy-Wolff Theorem. Also the step-lengths $s_{n}=d
(z_{n},z_{n+1})$ decrease to $s_{\infty }$, and 
$\phi $ is said to be of type I (or non-zero-step) if $s_{\infty }>0$;
$\phi $ is of type II (or zero-step) if $s_{\infty }=0$.
We will refine the type I and type II classes. 
By Julia's Lemma, $\ima z_{n}\uparrow L_{\infty }$. A map $\phi$ of type I
is said to be of {\bf type  Ia} (or non-zero-step/finite-height) if
$L_{\infty }<\infty $ and of  
{\bf type Ib} (or non-zero-step/infinite-height) if
$L_{\infty }=\infty$. Likewise for type II.
\begin{example}\label{ex:basic}
The map $\phi (z)=z+1$ is of type Ia, while $\phi (z)=z+i$ is of type IIb.
We give more examples in Section \ref{sec:examples}.
\end{example}
It is not clear {\em a priori} that this classification does not
depend on the choice of $i$ as starting point, however, that is indeed
the case, as we will see below, as a consequence of the next theorem.
\begin{theorem}[Pommerenke, \cite{pom} (3.17)]\label{thm:forward}
Let $\phi $ be an analytic self-map of $\Hh $ of parabolic type as in
(\ref{eq:assump}), and 
let $\{z_{n}=\phi _{n} (i)\}_{n=0}^{\infty }$ be a forward-iteration
sequence.  Then
\[
\frac{\ima z_{n+1}}{\ima z_{n}}\longrightarrow 1
\]
as $n$ tends to infinity.

Moreover, letting $z_{n}=u_{n}+iv_{n}$ and considering the
automorphisms of $\Hh $ given by $M_{n} (z)= (z-u_{n})/v_{n}$,
the normalized iterates 
$M_{n}\circ \phi _{n}$ converge uniformly on compact subsets of $\Hh $
to a function $\sigma $ which satisfies the functional equation
\[
\sigma \circ \phi =\sigma + b
\]
where 
\begin{equation}\label{eq:b}
b=\lim_{n\rightarrow \infty }\frac{u_{n+1}-u_{n}}{v_{n}}
\end{equation}
and $b\neq 0$ in the non-zero-step case, while $b=0$ in the zero-step case.
\end{theorem}
Since
\[
M_{n}\circ M_{n+1}^{-1}
(z)=z\frac{v_{n+1}}{v_{n}}+\frac{u_{n+1}-u_{n}}{v_{n}}\rightarrow z+b, 
\]
given and arbitrary point $z\in \Hh $,
\begin{eqnarray*}
d (\phi _{n} (z),\phi _{n+1} (z)) & = & d (M_{n}\circ\phi _{n}
(z),M_{n}\circ M_{n+1}^{-1}\circ  M_{n+1}\circ\phi _{n+1} (z) )\\
& \rightarrow &
d (\sigma (z),\sigma (z)+b). 
\end{eqnarray*}
Therefore the type I and type II classification is well-defined.
The fact that the type a and type b classification is well-defined
follows easily from Schwarz-Pick, because $d (\phi _{n} (z),\phi _{n}
(i))\leq d (z,i)$. 

\subsection{Parabolic backward iteration}\label{ssec:backward}

We now study parabolic maps $\phi $ as in (\ref{eq:assump}) which have
a BISBS tending to infinity.
The map $\phi (z)=z+i$ shows that such
sequences may not exist at all. We continue our classification by saying
that $\phi $ is {\bf type} $\mathbf \emptyset$ if it has no
BISBS. So, $\phi (z)=z+i$ is of type $\emptyset$. On the other
hand, if a BISBS exists $\{w_{n} \}_{n=0}^{\infty }$, then by Julia's
Lemma, $y_{n}:=\ima 
(w_{n})\downarrow \ell_{\infty }$. We say that a BISBS is of {\bf type}
$\mathbf 1$  (or non-zero-height) if  $\ell_{\infty }>0$, and
of {\bf type} $\mathbf 2$ (or zero-height) if
$\ell_{\infty }=0$. So, for instance, $\phi (z)=z+1$ has a BISBS of
type 1.
At first sight, one might think that type 2 never
arises. However, a simple example is given by the following map $\phi
(z)=\sqrt{z^{2}-1}$. Thinking of $\phi $ as the composition of three
simple operations, one checks that $\phi $ maps $\Hh $ into itself and
that it is of parabolic type IIb with Denjoy-Wolff point at infinity.
The sequence $w_{n}=\sqrt{n+i}$ is a backward-iteration sequence for
$\phi $. A calculation shows that the pseudo-hyperbolic steps $d_{n}$
stay bounded away from $1$. So $w_{n}$ is a BISBS and $\ima
w_{n}\asymp 1/\sqrt{n}$ tends to zero.
In this  example, although $y_{n}$ tends to zero, it does not do so
very fast, e.g., $\sum y_{n}=\infty $. This is a general fact.
\begin{lemma}\label{lem:sum}
Let $\phi $ be a parabolic map of $\Hh $ as above, and let $w_{n}$ be
a BISBS tending to infinity. If $y_{n}=\ima w_{n}$, then
\[
\sum _{n=0}^{\infty }y_{n}=\infty 
\]
\end{lemma}
\begin{proof}[Proof of Lemma \ref{lem:sum}]
Recall that $d_{n}=d (w_{n},w_{n+1})\rightarrow a<1$. Then 
\[
\frac{1}{d_{n}}=\left|1-\frac{2y_{n}i}{w_{n+1}-w_{n}}\right|\leq
1+2\frac{y_{n}}{|w_{n+1}-w_{n}|}. 
\]
So, 
\[
\liminf_{n\rightarrow \infty }\frac{y_{n}}{|w_{n+1}-w_{n}|}\geq
\frac{1-a}{2a}.  
\]
Letting $C_{0}= (1-a)/4a$, there is $n_{0}$ such that for $n\geq n_{0}$,
$y_{n}\geq C_{0}|w_{n+1}-w_{n}|$. Summing from $n_{0}$ to $N>n_{0}$,
\begin{eqnarray*}
\sum _{n=n_{0}}^{N }y_{n} & \geq & C_{0}\sum _{n=n_{0}}^{N
}|w_{n+1}-w_{n}|\geq C_{0} \left|\sum _{n=n_{0}}^{N
} (w_{n+1}-w_{n})\right|\\
& = & C_{0}|w_{N+1}-w_{n_{0}}|\rightarrow \infty 
\end{eqnarray*}
\end{proof}

A main question in this context is whether the ratios 
$y_{n+1}/y_{n}$ tend to $1$ or not. 
All we can infer so far, from Lemma \ref{lem:sum}
for example, is that 
\begin{equation}\label{eq:limsup}
\limsup_{n\rightarrow \infty }\frac{y_{n+1}}{y_{n}}=1.
\end{equation}

We will show that indeed $\lim _{n\rightarrow \infty }y_{n+1}/y_{n}=1$.
It will be useful, in the course of the proof, to consider the
hyperbolic steps of higher order. Namely, if $w_{n}$ is a BISBS, and
$k=1,2,3,\dots $, then
\begin{equation}\label{eq:highersteps}
d (w_{n},w_{n+k})\uparrow a_{k}.
\end{equation}
To see that the $a_{k}$ are also strictly less than one, consider the
hyperbolic distance 
\[
\rho =\log \frac{1+d}{1-d}
\]
Then $\rho (w_{n},w_{n+k})\uparrow \rho _{k}$, and by the triangle
inequality, $\rho _{k}\leq k\rho _{1}<\infty $. Thus, $a_{k}<1$.

\begin{remark}\label{rem:arg}
Note that $\Arg (w_{n})$ either tends to $0$ or to $\pi $.
In fact, $0<\ima w_{n}\leq \ima w_{0}$ and $w_{n}\rightarrow
\infty $. So, if $\Arg w_{n}$ were to oscillate between $0 $ and
$\pi  $, the sequence $w_{n}$ would accumulate everywhere on the
real axis, because of the bounded-steps condition, and this would
yield a contradiction. 
\end{remark}

The following theorem is a direct counter-part to Theorem \ref{thm:forward}.

\begin{theorem}\label{thm:ratios}
Let $\phi $ be an analytic self-map of $\Hh $ of parabolic type as in
(\ref{eq:assump}), and
let $\{w_{n}\}_{n=0}^{\infty }$ be a backward-iteration sequence with
bounded pseudo-hyperbolic steps $d_{n}=d (w_{n},w_{n+1})\uparrow
a<1$, which tends to infinity. Then
\[
\frac{\ima w_{n+1}}{\ima w_{n}}\longrightarrow 1
\]
as $n$ tends to infinity.
\end{theorem}

This is the key to proving the following conjugation result. 

\begin{theorem}\label{thm:conj}
Let $\phi $ be an analytic self-map of $\Hh $ of parabolic type as in (\ref{eq:assump}), and
let $\{w_{n}=x_{n}+iy_{n}\}_{n=0}^{\infty }$  be a backward-iteration
sequence with 
bounded pseudo-hyperbolic steps $d_{n}=d (w_{n},w_{n+1})\uparrow
a<1$, which tends to infinity. 
(Assume also WLOG that $\Arg w_{n}$ tends to $0$). 
Consider the automorphisms of $\Hh $ given by
$\tau_{n} (z)=x_{n}+y_{n}z$. Then the normalized iterates $\phi
_{n}\circ \tau_{n}$ converge uniformly on compact subsets of $\Hh $ to
an analytic self-map $\psi $ of $\Hh $ such that
\[
\psi (z-b_{0})=\phi \circ \psi (z)
\]
where 
\begin{equation}\label{eq:b0}
b_{0}=\frac{2a}{\sqrt{1-a^{2}}}=\lim _{n\rightarrow \infty }
\frac{x_{n+1}-x_{n}}{y_{n}}
\end{equation}
\end{theorem}

\begin{corollary}\label{cor:ratios}
With the hypothesis of Theorem \ref{thm:ratios} and Theorem
\ref{thm:conj}, letting $\zeta _{n}=i+nb_{0}$,
\begin{enumerate}
\item $\psi (\zeta _{n})=w_{n}$ and $\psi ^{\prime } (\zeta
_{n})/y_{n}\rightarrow 1$.
\item $\psi $ has non-tangential limit $\infty $ at $\infty $.
\item $\lim \phi ^{\prime } (w_{n})=1$.
\end{enumerate}
\end{corollary}

Theorem \ref{thm:ratios}, Theorem \ref{thm:conj} and Corollary
\ref{cor:ratios}  are proved in Section \ref{ssec:prf}.

\subsection{Uniqueness for BISBS of type 1}\label{ssec:uquests}

Let $w_{n}$ be a BISBS of non-zero-height as defined in Section
\ref{ssec:backward}. If $\tau _{n} (z)=x_{n}+zy_{n}$,
then Theorem \ref{thm:conj} 
says that $\phi _{n}\circ \tau _{n}$ tends to the conjugating map
$\psi$. So, given an arbitrary $z\in\Hh $, 
the sequence $\{\psi (z+nb_{0}) \}_{n=0}^{\infty }$ is a BISBS for
$\phi $ and is of non-zero-height as well, by Proposition
\ref{thm:conjtyp1} below and by Julia's Lemma. 
We show that every BISBS of
non-zero-height  occurs this way.
\begin{theorem}\label{thm:ufinheight}
Let $\phi $ be an analytic self-map of $\Hh $ of parabolic type, and
let $\{w_{n}=x_{n}+iy_{n}\}_{n=0}^{\infty }$  be a backward-iteration
sequence tending to the Denjoy-Wolff point, $\tau _{\phi }=\infty $, with 
bounded pseudo-hyperbolic steps $d_{n}=d (w_{n},w_{n+1})\uparrow
a<1$, and which is furthermore of non-zero-height, i.e.,
$y_{n}\downarrow \ell_{\infty }>0$.
Let $\tilde{w}_{n}$ be another backward-iteration
sequence tending to infinity with $d
(\tilde{w}_{n},\tilde{w}_{n+1})\uparrow \tilde{a}<1$, and
$\tilde{y}_{n}\downarrow \tilde{\ell}_{\infty }>0$. 
Assume further that $\Arg w_{n}$ and $\Arg \tilde{w}_{n}$ are both
tending to zero. Likewise let $\tau _{n}$, $\tilde{\tau }_{n}$, $\psi $, and
$\tilde{\psi }$, $b_{0}$, and $\tilde{b}_{0}$ be the corresponding
maps given by Theorem \ref{thm:conj}. Then
\[
\tilde{\psi } \left(\frac{\tilde{b}_{0}}{b_{0}}z+b_{1} \right)=\psi (z)
\]
for some $b_{1}\in\R $. So that,
$\tilde{w}_{n}=\psi (z_{1}+nb_{0})$
for $z_{1}= (b_{0}/\tilde{b}_{0}) (i-b_{1})$, and
\[
\lim _{n\rightarrow \infty }x_{n+1}-x_{n}=\ell_{\infty
}\frac{2a}{\sqrt{1-a^{2}}}=\tilde{\ell}_{\infty
}\frac{2\tilde{a}}{\sqrt{1-\tilde{a}^{2}}}= \lim _{n\rightarrow \infty
}\tilde{x}_{n+1}-\tilde{x}_{n}  
\]
\end{theorem}
\begin{question}\label{quest:twofinh}
It is not clear whether Theorem \ref{thm:ufinheight} holds in the
zero-height case, and whether 
zero-height BISBS can coexist with non-zero-height one. We conjecture
that at least for univalent maps $\phi $ this never happens. Theorem
\ref{thm:nonexist} below implies that this cannot happen for type
I parabolic maps.
\end{question}

Theorem \ref{thm:ufinheight} is proved in Section
\ref{ssec:unique}.

\subsection{More on BISBS of type 1}\label{ssec:conjtyp1}

\begin{proposition}\label{thm:conjtyp1}
Let $w_{n}=x_{n}+iy_{n}$ be a BISBS of non-zero-height, i.e. such that
$y_{n}\downarrow \ell_{\infty }>0$, and assume $x_{n}\rightarrow
+\infty$. Let
$\psi$ be the conjugating map obtained in Theorem \ref{thm:conj}. Then
\[
\frac{\psi (z)}{z}\longrightarrow \ell_{\infty } 
\]
as $z$ tends to infinity non-tangentially.
\end{proposition}
\begin{corollary}\label{cor:conjtyp1}
Let $\psi $ be as in Proposition \ref{thm:conjtyp1}, and consider the
quadrants $Q (R)=\{x>R,y>1 \}$. Then $\psi $ is one-to-one on $Q (R)$
for $R$ large enough, and this is true of $\phi $ as well. 
\end{corollary}

In the course of the proof of Theorem \ref{thm:ufinheight} above we
found useful to introduce the following notion of limit.
For $R>0$ and $\ep >0$ define the horizontal half-strips
\[
\St (R,\ep )=\{z=x+iy:\ x\geq R,\ \ep \leq y\leq 1/\ep  \}.
\]
\begin{definition}\label{def:laterallim}
Given a complex-valued function $f$ defined on $\Hh $, we say that $f$
has a {\sf lateral-limit at $+\infty$} if there is $a\in \C$ such that
\[
\lim _{\St (R,\ep )\ni z\rightarrow \infty }f (z)=a
\]
for every choice of $R>0$ and $0<\ep <1$. In this case we write,
\[
\latlim_{z\rightarrow +\infty }f (z)=a
\]
\end{definition}

\begin{proposition}\label{prop:philatlim}
Let $\phi $ be an analytic self-map of $\Hh $ of parabolic type as in
(\ref{eq:assump}), and 
let $\{w_{n}=x_{n}+iy_{n}\}_{n=0}^{\infty }$  be a backward-iteration
sequence with 
bounded pseudo-hyperbolic steps $d_{n}=d (w_{n},w_{n+1})\uparrow
a<1$, which tends to infinity,
and which is furthermore of non-zero-height, i.e., such that
$y_{n}\downarrow \ell_{\infty }>0$.
Assume also that $x_{n}$ tends to $+\infty $ and let $b_{0}$ be
defined as in (\ref{eq:b0}). Then,
\[
\latlim_{z\rightarrow +\infty } z-\phi (z)=b_{0}\ell_{\infty }
\]
\end{proposition}

The following question is therefore quite natural, also in view of Theorem
\ref{thm:finn} in the case of a BRFP. In that case, if $\phi $ fixes a
point $\zeta\neq \tau _{\phi } $ 
on the boundary of the unit disk, then the existence of
a finite angular derivative at $\zeta $ implies the existence of BISBS
tending to $\zeta $.
\begin{question}\label{quest:exist}
Suppose $\phi $ is an analytic self-map of $\Hh $, as
in (\ref{eq:assump}). Suppose moreover that
\[
\latlim _{z\rightarrow +\infty }z-\phi (z)=C_{0}>0
\]
Does this imply that there exists a BISBS of type 1 tending to $+\infty $?
\end{question}
Proposition \ref{thm:conjtyp1}, Corollary \ref{cor:conjtyp1}, and
Proposition \ref{prop:philatlim} are proved in Section \ref{ssec:typ1}.

\subsection{BISBS for type I parabolic maps}\label{ssec:nonexist}
We show that if a self-map $\phi $ of $\Hh $ is of
parabolic type I, then it cannot have a BISBS of
zero-height.
\begin{theorem}\label{thm:nonexist}
Let $\phi $ be an analytic self-map of $\Hh $, as
in (\ref{eq:assump}), which is of parabolic type I 
i.e. the forward iterates have non-zero step. 
Assume, also that $\phi $ admits a backward-iteration 
sequence $\{w_{n}\}_{n=0}^{\infty }$  
tending to infinity, with 
bounded pseudo-hyperbolic steps $d_{n}=d (w_{n},w_{n+1})\uparrow
a<1$. Then $w_{n}$ must be of non-zero-height, i.e.,
$\ima w_{n}\downarrow \ell_{\infty }>0$. 
\end{theorem}
\begin{corollary}\label{cor:typeIa}
A map $\phi $ as in (\ref{eq:assump}) of type Ia which has forward
iterates whose arguments tend to $0$ 
cannot have a BISBS whose argument also tends to $0$.
\end{corollary}
We don't know
whether Corollary \ref{cor:typeIa} holds for type Ib and type IIa maps. 
The fact that for type Ib maps the argument of the forward iterates
tends to either $0$ or $\pi $ is proved in Remark 1 of \cite{pom}.

Theorem \ref{thm:nonexist} and Corollary \ref{cor:typeIa} 
are proved in Section \ref{ssec:no}.

We now start with the proof of the various statements. In the last
section of the paper,
Section \ref{sec:examples}, we produce as many examples as we could
find of the different behaviors, however, some cases are missing.

\section{Convergence properties}\label{sec:convergence}
In this section we prove Theorem \ref{thm:converges}.
We begin with the proof of (1), which follows a standard line of
argument, however, the idea of using 
Theorem \ref{thm:cp} was suggested to us by F. Bracci.

\begin{proof}[Proof of Theorem \ref{thm:converges} (1)]
We treat first the elliptic case, $\tau_{\phi }\in\D$.
Since 
$d(w_n,\tau_{\phi }) \leq  d(w_{n+1}, \tau_{\phi })$,
$r=\lim _{n\rightarrow \infty }d(w_n,\tau_{\phi })$ exists.
If $r=0$, then $w_{n}=\tau _{\phi }$ for all $n$'s.
If $0<r<1$, either $\phi$ is an elliptic rotation or there exists
a constant $m<1$ such that $\phi $ maps the disk $\{z\in \D :d (z,\tau
_{\phi })< r \}$ into the disk
$\{z\in \D :d (z,\tau_{\phi })< mr \}$, and this yields a
contradiction. 
Therefore, aside for the trivial case when $\phi $ is a conjugate
rotation or $w_{n}$ is identically equal to $\tau _{\phi }$, we must
have $r=1$. Any subsequence, $w_{n_{k}}$ tending to a point $\zeta \in
\bd \D $ satisfies all the hypothesis of Corollary \ref{thm:jc}. 
Thus $\zeta $ is a BRFP for
$\phi $ with multiplier bounded by $A:=(1+a)/ (1-a)$.
However, the fact the $\{w_{n} \}$ has bounded steps implies that its
cluster set on $\bd \D $ must be connected, and hence is either a
point or an interval. To see this, connect $w_{n}$ to $w_{n+1}$ with a
straight segment to obtain a curve which has the same cluster set as
the sequence $\{w_{n} \}$ on $\bd \D$.
Such cluster set can't be an interval, because by Theorem
\ref{thm:cp}, the set of BRFPs
is at most countable (actually finite in the elliptic and hyperbolic
cases since the multipliers stay bounded). 
So there is a unique $\zeta \in\bd  \D $, such
that $w_{n}\rightarrow \zeta$ and
$\zeta $ is a BRFP for $\phi $, i.e. $1<\phi^\prime(\zeta)\leq A$.

We now treat the hyperbolic case. Without loss of generality, we
assume that $\tau_{\phi }=1$ and $0<\phi ^{\prime } (1):=c<1$.
Recall the Poisson kernel at $1$:
\begin{equation}\label{eq:poissonk}
P (w)=\frac{1-|w|^{2}}{|1-w|^{2}}
\end{equation}
Horocycles at $1$, i.e. disks interior to $\D$ and tangent to $\bd \D$
at $1$, can also be defined as level sets for $P$:
\[
H (t)=\{w: P (w) > 1/t\}.
\]
The existence of the angular derivative $\phi ^{\prime } (1)$ implies,
by Julia's Lemma (see \cite{sh} p. 63),  that
the horocycle $H(t)$ is mapped into the
horocycle $H (ct)$. Choose $t_{0}$ so small that $w_{0}$ does not belong to $H
(t_{0})$. Then,  $w_{n}$ is not in $H (c^{-n}t_{0})$, hence the
sequence $w_{n}$ cannot have cluster points in $\D$.
By the same arguments as in Case 1, we obtain that either $\{w_{n}
\}$ converges to a BRFP $\zeta$, or $\{w_{n} \}$
converges to $1$. Theorem \ref{thm:converges} (2) claims that the latter option
is impossible. 

Finally, in the parabolic case, assume $\tau_{\phi }=1$ and $\phi
^{\prime } (1):=c=1$. Note first that, $\{w_{n} \}$ cannot cluster in
$\D$. If so, there would be a subsequence $w_{n_{k}}$ tending to
$z_{0}\in\D$. But Schwarz-Pick implies that 
\[
d (\phi _{n_{k}} (z_{0}),w_{0})\leq d (z_{0},w_{n_{k}})\longrightarrow
0
\] 
and this contradicts the fact that $\phi _{n_{k}} (z_{0})$ tends to
$1$ by the Denjoy-Wolff Theorem.
Again, we conclude as before that  either $\{w_{n}
\}$ converges to a BRFP $\zeta$, or $\{w_{n} \}$
converges to $1$. This time however the latter option can occur. A
trivial example is the parabolic automorphism $z\mapsto z+1$ on the
upper half-plane (conjugated to the unit disk), but one can find other
examples as well, e.g. let $\sigma $ be a Riemann map of $\D$ onto
the set $\Omega=\{x+iy: y>\chi_{(-\infty ,0])} (x) \}$ which sends $1$
to $\infty $, and set $\phi(z) =\sigma ^{-1} (\sigma (z)+1)$. 
\end{proof}
We now prove
Theorem \ref{thm:converges} (2). In view of (1), 
it only remains to show that in the
hyperbolic case a BISBS cannot tend to $\tau _{\phi }$. This fact,
under the more stringent hypothesis of 
univalence on $\phi $, was already obtained by F. Bracci in \cite{br} by
different methods than ours. We thank Prof. Bracci for sharing is
preprint with us.

Before tackling the proof of Theorem \ref{thm:converges} (2), we want
to obtain an easy consequence of the bounded steps condition:
$d (w_{n},w_{n+1})\leq a<1$. 
The furthest $w_{n+1}$ can be from
the origin is $(|w_n|+a)/(1+a|w_n|)$. Therefore,
\[
1-|w_{n+1}|\geq\frac{1-a}{1+a}(1-|w_n|).
\]
On the other hand, the world's-greatest-identity tells us that
\begin{eqnarray*}
1-a^2 & \leq &
1-d(w_n,w_{n+1})^2=
\frac{(1-|w_n|^2)(1-|w_{n+1}|^2)}{|1-\overline{w_n}w_{n+1}|^2}\\
& \leq & 4a^2\frac{(1-|w_n|)(1-|w_{n+1}|)}{|w_n-w_{n+1}|^2}.
\end{eqnarray*}
Putting these two estimates together we obtain
\begin{equation}\label{eq:lb}
\frac{1-|w_{n+1}|}{|w_{n+1}-w_n|}\geq\frac{1-a}{2a}>0.
\end{equation}

\begin{lemma}\label{lem:ntapproach}
Let $\{w_{n} \}_{n=0}^{\infty }\subset \D$ be a sequence such that
$w_{n}\rightarrow \zeta \in \bd  \D$ as $n$ tends to infinity. Assume
that $w_{n}$      
satisfies (\ref{eq:lb}) and also
\begin{equation}\label{eq:angder}
\liminf_{n\rightarrow\infty}\frac{1-|w_n|}{1-|w_{n+1}|}\geq C>1.
\end{equation}
Then $w_{n}$ converges to $\zeta $ non-tangentially, i.e.,
\begin{equation}\label{eq:ntapproach}
\liminf_{n\rightarrow\infty}\frac{1-|w_n|}{|\zeta -w_n|}\geq \de>0
\end{equation}
\end{lemma}
\begin{proof}
Note that, 
\begin{eqnarray*}
|w_{n+1}|-|w_n| & = & (1-|w_n|)-(1-|w_{n+1}|)\\
& = &
(1-|w_{n+1}|)\left(\frac{1-|w_n|}{1-|w_{n+1}|}-1\right)
\end{eqnarray*}
Thus, by (\ref{eq:angder}),
\[
\liminf_{n\rightarrow\infty}\frac{|w_{n+1}|-|w_n|}{1-|w_{n+1}|}\geq C-1>0
\]
and by (\ref{eq:lb}),
\[
\liminf_{n\rightarrow\infty}\frac{|w_{n+1}|-|w_n|}{|w_{n+1}-w_n|}\geq (C-1)\frac{1-a}{2a}=\delta>0,
\]
i.e., there is $n_0$ such that for $k\geq n_0$:
\[|w_{k+1}|-|w_k|\geq (\delta/2)|w_{k+1}-w_k|\]
So by telescoping sums, for $m>n>n_0$,
\[
|w_m|-|w_n|\geq (\de/2)|w_m-w_n|
\]
Letting $m$ tend to infinity and then taking the $\liminf$ as $n$
tends to infinity, we find
that (\ref{eq:ntapproach}) holds.
\end{proof}

Now we are in position to prove (2) in Theorem \ref{thm:converges},
i.e., we need to show that in the hyperbolic case a backward-iteration
sequence with bounded hyperbolic steps cannot tend to the Denjoy-Wolff
point.
\begin{proof}[Proof of Theorem \ref{thm:converges} (2)]
Without loss of generality assume that $\tau_{\phi }=1$, so $0<c:=\phi
^{\prime } (1)<1$. Also assume that $w_{n}$ does tend to $1$, we will reach
a contradiction. 

First we rewrite (\ref{eq:lb}) as follows
\begin{equation}\label{eq:lb2}
\frac{1-a}{2a}\left|1-\frac{1-w_{n}}{1-w_{n+1}} \right|\leq
\frac{1-|w_{n+1}|}{|1-w_{n+1}|}
\end{equation}
Julia's Lemma tells us that since $w_{n}$ does not belong to the
horocycle $H (t)$ with $t=|1-w_{n}|^{2}/ (1-|w_{n}|^{2})$, $w_{n+1}$
cannot belong to the horocycle $H (t/c)$, i.e.,
\begin{equation}\label{eq:jl}
\frac{1-|w_{n+1}|^{2}}{|1-w_{n+1}|^{2}}\leq
c\frac{1-|w_{n}|^{2}}{|1-w_{n}|^{2}} 
\end{equation}
Iterating this estimate we find that
\begin{equation}\label{eq:tlim}
\lim_{n\rightarrow \infty }\frac{1-|w_{n}|^{2}}{|1-w_{n}|^{2}}=0
\end{equation}
Applying this to (\ref{eq:lb2}) we obtain that
\[
\lim _{n\rightarrow \infty }\frac{1-w_{n}}{1-w_{n+1}}=1
\]
Going back to (\ref{eq:jl}),
\[
\frac{1-|w_{n}|^{2}}{1-|w_{n+1}|^{2}}\geq
\frac{1}{c}\left|\frac{1-w_{n}}{1-w_{n+1}}\right| 
\]
Therefore,
\[
\liminf_{n\rightarrow \infty }\frac{1-|w_{n}|}{1-|w_{n+1}|}=
\liminf_{n\rightarrow \infty }\frac{1-|w_{n}|^{2}}{1-|w_{n+1}|^{2}}\geq\frac{1}{c}>1.
\]
Now we can apply Lemma \ref{lem:ntapproach}, with $C$ replaced by
$1/c$, and deduce that (\ref{eq:ntapproach}) holds with $\zeta =1$.
This contradicts (\ref{eq:tlim}).
\end{proof}

\begin{proof}[Proof of Theorem \ref{thm:converges} (3)]
Suppose now that $w_{n}$ is a BISBS
converging to a BRFP $\zeta \in\bd \D$. Without loss of generality
$\zeta =1$.  
Let $A:=\phi ^{\prime } (1)>1$.
By Theorem \ref{thm:juca}, 
\[
A=\liminf_{z\rightarrow \zeta}\frac{1-|\phi (z)|}{1-|z|}.
\]
In particular, $w_{n}$ satisfies (\ref{eq:angder}) with
$C=A$. Therefore, by Lemma \ref{lem:ntapproach}, we find that
(\ref{eq:ntapproach}) holds.

Now let $\psi $ be the conjugation produced by Theorem \ref{thm:finn}. 
In Lemma 5.1 of \cite{finn} we show the existence of a simply
connected region $\Omega\subset \Hh$ with an inner tangent at $0$ with
respect to $\Hh$, such
that $\psi $ is one-to-one on $\Omega$, and so that $\psi (\Om)$ has an
inner tangent at $1$ with respect to $\D$. 
So, eventually, $w_{n}\in\psi  (\Om)$. Let $f$ be
the inverse of $\psi_{\mid\Om} $. Then 
$b_{n}:=f(w_{n})=A^{-n}b_{0}$ for some $b_{0}\in\Hh $.
By (\ref{eq:semiconf}), this proves (\ref{eq:angappr}).
\begin{remark}\label{rem:typo}
We found that the proof of Lemma 5.1 of \cite{finn} has a typo, i.e.,
$\theta_{n}\uparrow \pi /2$ instead of $\theta _{n}\downarrow 0$.
\end{remark}

\end{proof}

\begin{proof}[Proof of Theorem \ref{thm:converges} (4)]
Assume now that $\tau_{\phi }=1$ and $c:=\phi ^{\prime }
(1)=1$. Assume also that $\{w_{n} \}_{n=0}^{\infty }$ is a BISBS
coverging to $1$. By Julia's Lemma, $P (w_{n+1})\leq P (w_{n})$, where
$P$ is defined in (\ref{eq:poissonk}). This shows that $w_{n}$ tends
to $1$ tangentially.
\end{proof}

\section{The parabolic case}\label{sec:parabolic}

\subsection{Proof of Theorem \ref{thm:ratios}, Theorem \ref{thm:conj}  and 
Corollary \ref{cor:ratios}}\label{ssec:prf}  
We write
$w_{n}=x_{n}+iy_{n}$ for simplicity. Note that if $y_{n}\downarrow
\ell_{\infty }>0$, the fact that $y_{n+1}/y_{n}\rightarrow 1$ is clear
in this case. However, we will treat both cases together.
Recall that $d_{n}=d (w_{n},w_{n+1})\uparrow a<1$ and that $\tau _{n}
(z)=x_{n}+zy_{n}$. 
Notice that
\begin{equation}\label{eq:tausteps}
\tau _{n}^{-1}\circ \tau_{n+1}
(z)=\frac{x_{n+1}-x_{n}}{y_{n}}+\frac{y_{n+1}}{y_{n}}z 
\end{equation}
is another automorphism of $\Hh $, and $d_{n}=d (i,\tau _{n}^{-1}\circ
\tau _{n+1} (i))$. 
Hence, we obtain after some
manipulation: 
\begin{equation}\label{eq:bddstep}
\left(\frac{x_{n+1}-x_{n}}{y_{n}}
\right)^{2}=-\left(\frac{y_{n+1}}{y_{n}}
\right)^{2}+2\frac{1+d_{n}^{2}}{1-d_{n}^{2}} \left(\frac{y_{n+1}}{y_{n}}
\right)-1.
\end{equation}
The right hand-side being a concave-down quadratic polynomial in
$y_{n+1}/y_{n}$ implies that 
\begin{equation}\label{eq:xsteps}
\frac{|x_{n+1}-x_{n}|}{y_{n}}\leq \frac{2d_{n}}{1-d_{n}^{2}}\leq
\frac{2a}{1-a^{2}}   
\end{equation}
and, the left hand-side being positive yields
\begin{equation}\label{eq:ysteps}
\frac{1-a}{1+a}\leq \frac{(1-d_{n})^{2}}{1-d_{n}^{2}}=
\frac{1+d_{n}^{2}}{1-d_{n}^{2}}-\frac{2d_{n}}{1-d_{n}^{2}}\leq 
\frac{y_{n+1}}{y_{n}}\leq 1.
\end{equation}
\begin{claim}\label{cl:identity}
The following holds:
\[
f_{n} (z):=\frac{\phi \circ \tau_{n+1} (z)-w_{n}}{\phi \circ \tau_{n+1}
(z)-\overline {w_{n}}} \cdot\frac{z+i}{z-i}\longrightarrow 1  
\]
as $n$ tends to infinity, uniformly on compact subsets of $\Hh $.
\end{claim}
\begin{proof}[Proof of Claim \ref{cl:identity}]
Since $\phi (\tau _{n+1} (i))=\phi (w_{n+1})= w_{n}$, $f_{n} (z)$ is a
well-defined analytic function for $z\in\Hh $, and by Schwarz-Pick,
\begin{equation}\label{eq:fnless1}
|f_{n} (z)|=\frac{d (\phi \circ \tau _{n+1} (z),\phi \circ \tau _{n+1}
(i))}{d (z,i)}\leq 1 
\end{equation}
for all $z\in \Hh $. 

We first show that $|f_{n} (z)|$ tends to $1$ uniformly on compact
subsets of $\Hh $.
Consider a subsequence $f_{n_{j}}$.
By normal families, we can extract a subsequence $f_{N}$ tending to $f$.
By (\ref{eq:xsteps}) and (\ref{eq:ysteps}), we can extract a
subsequence $M$ of $N$ 
so that $y_{M+2}/y_{M+1}\rightarrow c\leq 1$ and
$(x_{M+2}-x_{M+1})/y_{M+1}\rightarrow b\geq 0$ ($b$ can be
chosen to be positive because we can
assume without loss of generality that $\Arg w_{n}$ tends to $0$, in
view of Remark \ref{rem:arg}). Thus,
$\tau _{M+1}^{-1}\circ \tau_{M+2}$ converges uniformly
on compact subsets of $\C$ to the automorphism $S (z)=cz+b$, and
by (\ref{eq:bddstep}),
\begin{equation}\label{eq:etacoeff}
b^{2}=-c^{2}+2\frac{1+a^{2}}{1-a^{2}}c-1.
\end{equation}
By Schwarz-Pick again,
\begin{equation}\label{eq:zerolim}
d (\phi \circ \tau_{M+1}\circ S(i),\phi \circ \tau_{M+2}
(i))\leq 
d (S (i),\tau _{M+1}^{-1}\circ \tau _{M+2} (i))\rightarrow 0
\end{equation}
as $M$ tends to infinity. Using the triangle inequality for the
hyperbolic distance, and then transferring it back to the
pseudo-hyperbolic distance (this works because of (\ref{eq:zerolim})), we obtain
\begin{eqnarray*}
|f_{M} (S (i))| & = & \frac{d (\phi \circ \tau_{M+1}\circ
S(i),\phi \circ \tau_{M+1} (i))}{d (S (i),i)}\\
& \geq & \frac{d (\phi \circ \tau_{M+2} (i),\phi \circ \tau_{M+1}
(i))}{d (S (i),i)}- o (1)\\
& = & \frac{d (w_{M},w_{M+1})}{d (S (i),i)}-o (1) 
\end{eqnarray*}
On the other hand, using (\ref{eq:etacoeff}),
\begin{equation}\label{eq:zetasteps}
d (S (i),i)^{2} = \frac{b^{2}+ (c-1)^{2}}{b^{2}+
(c+1)^{2}}=a^{2}
\end{equation}
Thus, since $d_{n}$ tends to $a$,
\[
\lim _{M\rightarrow \infty}|f_{M} (S (i))|=1.
\]
Thus $f$ is a constant of modulus one, and 
\[
\lim _{n\rightarrow \infty}|f_{n} (z)|=1
\]
uniformly on compact subsets of $\Hh $.

On the other hand,
\begin{equation}\label{eq:fnati}
f_{n} (i)= (\phi\circ \tau _{n+1}) ^{\prime } (i)\cdot
\frac{2i}{w_{n}-\overline {w_{n}}}=\phi ^{\prime } (w_{n+1})\cdot
\frac{y_{n+1}}{y_{n}} 
\end{equation}
So
\begin{equation}\label{eq:limmodfi}
\lim _{n\rightarrow \infty }\frac{y_{n+1}}{y_{n}}|\phi ^{\prime } (w_{n+1})|=1.
\end{equation}
Write $\phi (z)=z+p (z)$ with $\ima p (z)\geq 0$ and
$p (z)\rightarrow 0$ as $z\rightarrow \infty $ non-tangentially. Note that,
\begin{equation}\label{eq:limimpy}
\frac{\ima p (w_{n+1})}{y_{n+1}}=\frac{\ima
(w_{n}-w_{n+1})}{y_{n+1}}=\frac{y_{n}}{y_{n+1}}-1
\end{equation}
Also Schwarz-Pick applied to $p$ yields,
\[
\frac{y_{n+1}}{y_{n}}|p^{\prime } (w_{n+1})|\leq\frac{y_{n+1}}{y_{n}}
\frac{\ima p (w_{n+1})}{\ima w_{n+1}} =1-\frac{y_{n+1}}{y_{n}}
\]
for all $n $.
Therefore,
\begin{eqnarray*}
\left(\frac{y_{n+1}}{y_{n}} \right)^{2}|\phi ^{\prime } (w_{n+1})|^{2}
& = &\left(\frac{y_{n+1}}{y_{n}} \right)^{2}\left( 1+|p^{\prime }
(w_{n+1})|^{2}+2\rea 
p^{\prime } (w_{n+1}) \right)\\  
& \leq  & 2\frac{y_{n+1}}{y_{n}}\left(\frac{y_{n+1}}{y_{n}}\rea \phi
^{\prime } (w_{n+1})-1 \right)+1 
\end{eqnarray*}
Rearraging this inequality,
using the fact that $(y_{n+1}/y_{n})\rea \phi ^{\prime } (w_{n+1})$ is
less than $1$ (by (\ref{eq:fnless1}) applied to $z=i$), and the fact that
$y_{n+1}/y_{n}$ is greater than $(1-a)/ (1+a)$, see
(\ref{eq:ysteps}), we find that
\[
0\leq 2\frac{1-a}{1+a}\left(1-\frac{y_{n+1}}{y_{n}}\rea \phi ^{\prime
} (w_{n+1}) \right)\leq 1-\left(\frac{y_{n+1}}{y_{n}} \right)^{2}|\phi
^{\prime } (w_{n+1})|^{2} 
\]
By (\ref{eq:limmodfi}), we obtain
\[
\lim _{n\rightarrow \infty }\frac{y_{n+1}}{y_{n}}\rea \phi ^{\prime
} (w_{n+1})  =1.
\]
Therefore, $f_{n} (i)\rightarrow 1$, and  thus $f_{n} (z)\rightarrow
1$ on compact subsets of  $\Hh $. So
Claim \ref{cl:identity} is proved.
\end{proof}
Renormalize the iterates of $\phi $ by writing $\psi_{n} (z) =\phi
_{n}\circ \tau_{n} (z)$. Note that $\psi_{n} (i)=w_{0}$ and $\ima
\psi_{n}>0$. By Claim \ref{cl:identity},
\[
f_{n} (z)=\frac{\tau _{n}^{-1}\circ \phi \circ \tau _{n+1} (z)-i}{\tau
_{n}^{-1}\circ \phi \circ \tau _{n+1} (z)+i}\cdot \frac{z+i}{z-i}\rightarrow 1 
\]
which implies 
\begin{equation}\label{eq:tauphitau}
\lim _{n\rightarrow \infty}\tau _{n}^{-1}\circ \phi \circ \tau _{n+1} (z)=z.
\end{equation}
So, by Schwarz-Pick applied to $\phi _{n}$ and conformal invariance,
\begin{equation}\label{eq:dsteps0}
d (\psi_{n+1} (z),\psi_{n} (z)) \leq d (\tau _{n}^{-1}\circ \phi \circ \tau
_{n+1} (z),z) \rightarrow 0
\end{equation}
as $n$ tends to infinity.

Let $y_{N}$ be a
subsequence of the $y_{n}$ such that the ratios $y_{N+1}/y_{N}$
converge to a constant $c$. At the moment, we only know that 
\[
0<\frac{1-a}{1+a}\leq c\leq 1.
\]
Recall that $\tau _{N}^{-1}\circ \tau _{N+1}$ tends to an
automorphism $S (z)=cz+b$ of $\Hh $, where $b$ is determined by
(\ref{eq:etacoeff}). 
By normal families we can assume, passing to a subsequence, that the
corresponding sequence of normalized iterates
$\psi_{N}$ tends to an analytic function $\psi $ uniformly on compact
subsets of $\Hh $. Then, $\psi (i)=w_{0}$ and $\ima \psi (z)>0$ for
all $z\in\Hh $. 
Note that 
\begin{eqnarray*}
d (\psi _{N} (S (i)),w_{1}) & = & d (\phi _{N}\circ \tau _{N}\circ S
(i),\phi _{N}\circ \tau _{N+1} (i))\\ 
& \leq & d (S (i),\tau _{N}^{-1}\circ \tau _{N+1} (i))\rightarrow 0
\end{eqnarray*}
as $N$ tends to infinity. Thus, $\psi (S (i))=w_{1}\neq w_{0}$, and
$\psi $ is not constant.
By (\ref{eq:dsteps0}), the sequence $\psi _{N+1}$ tends to the same
function $\psi $, and since $\psi _{N+1}=\phi \circ \psi _{N}\circ
(\tau _{N}^{-1}\circ \tau _{N+1})$, we find that $\psi $ must satisfy
the functional equation 
\begin{equation}\label{eq:functeq}
\psi =\phi \circ \psi \circ S.
\end{equation}
We now consider the sequence $\psi _{N+2}$, which also tends to $\psi$.
Again $\psi _{N+2}=\phi \circ \psi _{N+1}\circ
(\tau _{N+1}^{-1}\circ \tau _{N+2})$. By the same arguments as before,
see (\ref{eq:xsteps}) (\ref{eq:ysteps}) and (\ref{eq:tausteps}), given
a subsequence of $N$ we can extract another subsequence, which we call
$M$, so that $\tau _{M+1}^{-1}\circ \tau _{M+2}$ converges to an
automorphism $\tilde{S} (z)=\tilde{c}z+\tilde{b}$. Then, $\tilde{S}$
satisfies
\[
\phi \circ \psi \circ \tilde{S}=\psi =\phi \circ \psi \circ S.
\]
Since $\phi \circ \psi $ is non-constant we can invert it locally, and
since $S$ and $\tilde{S}$ are linear, they must coincide.
In particular, $\tau _{N+1}^{-1}\circ \tau _{N+2}$ must converge to
$S$ as well. A similar argument yields, for every $k=0,1,2,3\dots $, 
that $\tau _{N+k}^{-1}\circ \tau _{N+k+1}$ tends to $S$. Therefore,
$\tau _{N}^{-1}\circ \tau _{N+k}$ must tend to $S_{k}=S\circ \cdots
\circ S$, $k$ times. Write $\zeta _{k}=S_{k} (i)$. Then
\[
d (\psi_{N} (\zeta _{k}),w_{k})=d (\phi _{N}\circ \tau _{N} (\zeta
_{K}),\phi _{N}\circ \tau _{N+k} (i))\leq  d (\zeta _{k},\tau
_{N}^{-1}\circ \tau _{N+k} (i))\rightarrow 0
\]
that is to say 
\begin{equation}\label{eq:psizetak}
\psi (\zeta _{k})=w_{k}.
\end{equation}

We now consider the functions
\[
g_{n} (z)=\frac{\psi \circ S_{n} (z)-w_{n}}{\psi \circ S_{n}
(z)-\overline {w_{n}}} \cdot \frac{z+i}{z-i} 
\]
Since $\psi \circ S_{n} (i)=w_{n}$, the $g_{n}$ are analytic on $\Hh $.
Moreover, by Schwarz-Pick,
\[
|g_{n} (z)|=\frac{d (\psi \circ S_{n} (z),\psi \circ S_{n} (i))}{d
(z,i)}\leq 1 
\]
Note that by (\ref{eq:zetasteps})
\[
|g_{n} (\zeta _{1})|=\frac{d (w_{n+1},w_{n})}{d
(\zeta _{1},i)}\rightarrow \frac{a}{d(\zeta _{1},i)}=1.
\]
Therefore, any subsequence of $|g_{n} (z)|$ has a subsequence
converging to a constant of modulus one, i.e.,
\[
|g_{n} (z)|\longrightarrow 1
\]
as $n$ tends to infinity, for all $z\in\Hh $.
Evaluating $|g_{n}|$ at $\zeta _{k}$, we obtain
\[
\frac{d (w_{n+k},w_{n})}{d (\zeta _{k},i)}\longrightarrow 1
\]
as $n$ tends to infinity. 
Thus, the number $a_{k}$ introduced in (\ref{eq:highersteps}) satisfy
\[
a_{k}=d (\zeta _{k},i)
\]
Recall that $\zeta _{k}=S_{k} (i)$ and $S (z)=cz+b$. 

Now assume that $c<1$, then 
\[
\zeta _{k}=c^{k}i+\frac{1-c^{k}}{1-c}b\longrightarrow\frac{b}{1-c} 
\]
Therefore, for $k$ large,
\[
1-a_{k}^{2}=\frac{4c^{k}}{|\zeta _{k}+i|^{2}}\leq 4c^{k}
\]
On the other hand, if $c=1$, then $\zeta _{k}=i+kb_{0}$, where by
(\ref{eq:etacoeff}), 
\begin{equation}\label{eq:bnot}
b_{0}=\frac{2a}{\sqrt{1-a^{2}}}.
\end{equation}
Thus, for $k$ large,
\[
1-a_{k}^{2}=\frac{4}{|\zeta
_{k}+1|^{2}}=\frac{4}{4+b_{0}^{2}k^{2}}\geq \frac{C}{k^{2}} 
\]
for some constant $C>0$.

These two asymptotic behaviors of the numbers $a_{k}$, as $k$ tends to
infinity, show that the sequence of ratios $y_{n+1}/y_{n}$ must
converge, and  (\ref{eq:limsup}) implies that the limit must be one.
So Theorem \ref{thm:ratios} is proved.

Therefore we have 
\[
\tau _{n}^{-1}\circ \tau _{n+1}\rightarrow z+b_{0}:=S (z)
\]
where $b_{0}$ is given in (\ref{eq:bnot}). Moreover, letting $\psi $
be a normal limit of the $\psi _{n}=\phi _{n}\circ \tau _{n}$, we have
have $\psi =\phi _{n}\circ \psi \circ S_{n}$. Also, letting 
\begin{equation}\label{eq:hn}
h_{n} (z)=\tau _{n}^{-1}\circ \psi \circ S_{n} (z)
\end{equation}
we find that
\[
h_{n} (\zeta _{k})=\tau _{n}^{-1}\circ \tau_{n+k} (i)\rightarrow \zeta _{k}
\]
So $h_{n}$ must tend to the identity, and 
by the same argument as in \cite{finn},
\begin{eqnarray*}
d (\phi _{n}\circ \tau _{n},\psi ) & = & d(\phi _{n}\circ \tau
_{n},\psi\circ S_{n}^{-1}\circ S_{n} )\\
& = &  d(\phi _{n}\circ \tau
_{n},\phi _{n}\circ \psi\circ S_{n} )\\
& \leq  &   d(\tau
_{n},\psi\circ S_{n} )\rightarrow 0.
\end{eqnarray*}
This proves Theorem \ref{thm:conj}.

Next we show Corollary \ref{cor:ratios}. The fact that $\psi (\zeta
_{n})=w_{n}$ is (\ref{eq:psizetak}). Also, by differentiating $h_{n}
(z)\rightarrow  z$, see (\ref{eq:hn}), we find that
\[
\frac{\psi^{\prime } (\zeta _{n})}{y_{n}}\rightarrow 1
\]
as $n$ tends to infinity. This is Corollary \ref{cor:ratios} (1). In
particular, $\psi ^{\prime } (\zeta _{n})\neq 0$ for $n$ large.

Now consider the half-line $\ga =\cup _{n=0}^{\infty }[\zeta
_{n},\zeta _{n+1}]$. Then $\psi $ tends to infinity along $\ga $
because $\psi (\zeta _{n})$ does and for $\zeta \in [\zeta
_{n},\zeta _{n+1}]$,
\[
d (\psi (\zeta ),\psi (\zeta _{n}))\leq Const.
\] 
So by Lindel\"{o}ff's Theorem, see \cite{pom2} Cor. 2.17 (i), $\psi$
has non-tangential limit $\infty $ at $\infty $.

Finally, Corollary \ref{cor:ratios} (3) follows from differentiating
(\ref{eq:tauphitau}).

\subsection{Uniqueness in the non-zero-height case}\label{ssec:unique}
Here we show Theorem \ref{thm:ufinheight} and so we refer to the
assumptions made in the statement.  
The next Lemma
amounts to say that $\psi (z)$ is asymptotic to the map
$\ell_{\infty }z$ as $z$ tends to $\infty $ 
``laterally''. For the definition of $\latlim _{z\rightarrow +\infty
}$ see Definition \ref{def:laterallim}.
\begin{lemma}\label{lem:latlim}
With $\psi $ as in Theorem \ref{thm:ufinheight}, the following holds,
\begin{equation}\label{eq:latlimim}
\latlim_{z\rightarrow +\infty }\frac{\ima \psi (z)}{\ima z}=\ell_{\infty }
\end{equation}
and
\begin{equation}\label{eq:latlimrea}
\latlim_{z\rightarrow +\infty }\frac{\rea \psi (z)}{\rea z}=\ell_{\infty }
\end{equation}
Moreover, we also have
\begin{equation}\label{eq:latlimder}
\latlim_{z\rightarrow +\infty }\psi ^{\prime} (z)=\ell_{\infty }.
\end{equation}
\end{lemma}
\begin{proof}
Recall that 
\[
y_{n}\rightarrow \ell_{\infty }\qquad \mbox { and }\qquad
x_{n+1}-x_{n}\rightarrow  b_{0}\ell_{\infty }
\]
So given $\varepsilon >0$ there is $n_{0}$ such that for $n\geq n_{0}$
\[
 b_{0}\ell_{\infty }-\varepsilon \leq x_{n+1}-x_{n}\leq
b_{0}\ell_{\infty }+\varepsilon  
\]
So 
\[
x_{n_{0}}+ (n-n_{0}) (b_{0}\ell_{\infty }-\varepsilon)\leq x_{n}\leq
x_{n_{0}}+ (n-n_{0}) (b_{0}\ell_{\infty }+\varepsilon) 
\]
Hence
\begin{equation}\label{eq:xn}
\lim _{n\rightarrow \infty }\frac{x_{n}}{nb_{0}}=\ell_{\infty }.
\end{equation}
Now recall that $h_{n} (z)=\tau_{n}^{-1}\circ \psi \circ S_{n} (z)$
tends to $z$ 
uniformly on compact subsets of $\Hh $, see (\ref{eq:hn}). Therefore
\[
\lim _{n\rightarrow \infty }|\psi (z+nb_{0})- (x_{n}+zy_{n})|=0
\]   
In particular,
\[
\lim _{n\rightarrow \infty }|\ima \psi (z+nb_{0})- \ell_{\infty }\ima z|=0
\]
uniformly for $z\in \{0\leq x\leq b_{0}, \ep \leq y\leq 1/\ep
\}$. This implies (\ref{eq:latlimim}).
And also,
\[
\lim _{n\rightarrow \infty }|\rea \psi (z+nb_{0})- x_{n}-\rea
z\ell_{\infty }|=0 
\]
which if we divide by $\rea z+nb_{0}$ still tends to zero. So, by
(\ref{eq:xn}), 
\[
\lim _{n\rightarrow \infty }\frac{\rea \psi (z+nb_{0})}{\rea
z+nb_{0}}=\ell_{\infty }. 
\]
This implies (\ref{eq:latlimrea}).
Finally, by differentiating $\tau_{n}^{-1}\circ \psi \circ S_{n} (z)$
we find that $\psi ^{\prime } (z+nb_{0})/y_{n}$ tends to $1$ uniformly
on compact subsets of $\Hh $, which implies (\ref{eq:latlimder}).
\end{proof}

As a corollary we obtain the following result which is similar to
Lemma 5.1 of \cite{finn}.
First we need another definition.
\begin{definition}\label{def:innerlat}
Given a simply connected region $\Om$ in $\Hh $, we say that $\Om$ has
an {\sf inner-lateral tangent} at $+\infty $, if for every choice of $0<\ep<1 $
there is $R>0$ large enough so that the half-strip
$\St (R,\ep )$ is contained in $\Om$.
\end{definition}
\begin{corollary}\label{cor:latlim}
Let $\psi$ be the conjugation obtained in Theorem \ref{thm:conj} for a
BISBS of type 1 whose argument tends to $0$. Then,
there is a convex region $\Om\subset \Hh $ with an inner-lateral
tangent at $+\infty $, such that $\psi $ is one-to-one on $\Om$ and
$\psi (\Om)$ also has an inner-lateral tangent at $+\infty$.
\end{corollary}
\begin{proof}
We build $\Om$ by induction on $k=2,3,\dots $. For $k=2$, choose
$R_{2}>0$ so that the three limits in Lemma \ref{lem:latlim} are within
$\ell_{\infty }/2$ of $\ell_{\infty }$ for $z\in\St (R_{2},1/2)$, and
set $p_{2}^{+}=p_{2}^{-}=R_{2}+i$. 
For $k$ arbitrary, consider the half-strips $\St ( R,1/k)$
with $R>R_{k-1}$ and let $p_{k}^{+}$ and $p_{k}^{-}$ be the two
intersection points between the boundary of  $\St ( R,1/k)$ and the
boundary of $\St ( R_{k-1},1/ (k-1))$ (The plus being assigned to the
point which is above the other one). Again choose $R$ large enough so
that three limits in Lemma \ref{lem:latlim} are within
$\ell_{\infty }/2$ of $\ell_{\infty }$ for $z\in\St (
R,1/k)$. Moreover, choose $R$ so large that the slope of the 
interval $[p_{k-1}^{+},p_{k}^{+}]$ is smaller than the slope of the
previous interval $[p_{k-2}^{+},p_{k-1}^{+}]$, likewise, so that the
slope of $[p_{k-1}^{-},p_{k}^{-}]$ is larger than that of
$[p_{k-2}^{-},p_{k-1}^{-}]$. 

Now let $\Om$ be the convex hull of the points $p_{k}^{\pm}$, $k=2,3,\dots $.
The univalence of $\psi $ on $\Om$ follows from Proposition 1.10 of
\cite{pom2}, because $\rea \psi ^{\prime }>0$ there and $\Om$ is
convex. Also it follows by construction that $\psi (\Om)$ contains
half-strips $\St (R,\ep )$ of arbitrarily large height.
So Corollary \ref{cor:latlim} is proved.
\end{proof}

\begin{proof}[Proof of Theorem \ref{thm:ufinheight}]
We proceed as in the uniqueness part of Theorem 1.2 of \cite{finn}.
Suppose $\psi$ and $\tilde{\psi }$ are the two conjugations.
Let $\Om$ and $\tilde{\Om}$ be the corresponding sets given by
Corollary \ref{cor:latlim}. Given a small hyperbolic disk centered at $i$,
$\De\subset \overline {\De }\subset \Hh $, there is $N$ such
that for all $n\geq  N$, $\De +nb_{0}\in \Om$, and $\psi
(\De +nb_{0})\subset  \tilde{\psi }  (\tilde{\Om})$. So, 
for $z\in \De $, $z+nb_{0}\in \Om$ and $\psi (z+nb_{0})\in
\tilde{\psi } (\tilde{\Om})$. Let $f$ denote the inverse of $\tilde{\psi }$ on
$\tilde{\psi } (\tilde{\Om})$ which maps $\tilde{\psi } (\tilde{\Om})$
back to $\tilde{\Om}$.  Then,
\[
\beta (z)=f (\psi (z+nb_{0}))-n\tilde{b}_{0}
\]
is analytic on $\De $, one-to-one there, 
and does not depend on $n\geq N$. In fact, $f
(\psi(z+nb_{0}))$ and $f (\psi (z+nb_{0}+b_{0}))$ are
mapped by $\tilde{\psi }$ to two points $z_{1}$ and $z_{2}$ such that
$\phi (z_{2})=z_{1}$. So,
\[
f (\psi (z+nb_{0}+b_{0}))-f(\psi (z+nb_{0}))=\tilde{b}_{0}.
\]
Since the radius of $\De $ was arbitrary, this implies that $\beta $
is well-defined and one-to-one on all of $\Hh $. By interchanging the
role of $\psi$ and $\tilde{\psi }$ we find another function $\ga $
analytic and one-to-one on all of $\Hh $ which is locally the inverse
of $\beta$, and thus also globally. So $\beta $ is an automorphism of $\Hh $.
Since $\beta =\tilde{S}_{n}^{-1}\circ \tilde{\psi }^{-1}\circ \psi\circ
S_{n}$, we have $\beta \circ S=\tilde{S}\circ \beta $, i.e., $\beta
(z+nb_{0})=\beta (z)+n\tilde{b}_{0}$. So $\beta $ fixes infinity and no other
point, i.e., $\beta (z)=c_{1}z+b_{1}$ with $b_{1}, c_{1}\in\R $.
In particular, $c_{1} (z+b_{0})+b_{1}=c_{1}z+b_{1}+\tilde{b}_{0}$, so
$c_{1}=\tilde{b}_{0}/b_{0}$. It follows from the definition of $\beta
$ that $\tilde{\psi }\circ \beta =\psi $. Hence,
\[
\tilde{\psi } \left(\frac{\tilde{b}_{0}}{b_{0}}z+b_{1} \right)=\psi (z)
\]
With this identity it is easy to check that 
\[
\tilde{w}_{n}=\psi (z_{1}+nb_{0})
\]
where $z_{1}= (b_{0}/\tilde{b}_{0}) (i-b_{1})$.

Finally, using (\ref{eq:latlimder}) we find that
\[
\ell_{\infty }b_{0}=\tilde{\ell}_{\infty }\tilde{b}_{0}
\]
\end{proof}

\subsection{Further properties for BISBS of type 1}\label{ssec:typ1}

\begin{proof}[Proof of Proposition \ref{thm:conjtyp1}]
Changing variables from $\Hh $ to $\D$ via the map $\alpha (z)= (z-i)/
(z+i)$, we find that $\Psi=\alpha \circ \psi \circ \alpha ^{-1}$ is an
analytic self-map of $\D$ with non-tangential limit $1$ at $1$, and a
calculation shows that
\[
\frac{1-|\Psi (\alpha (z) )|^{2}}{1-|\alpha (z) |^{2}}=\frac{\ima \psi
(z)}{\ima z} \frac{|z+i|^{2}}{|\psi (z)+i |^{2}}
\]
So when $z=\zeta _{n}$,
\[
\frac{1-|\Psi (\alpha (\zeta _{n}) )|^{2}}{1-|\alpha (\zeta _{n})
|^{2}}=y_{n}\frac{4+n^{2}b_{0}^{2}}{x_{n}^{2}+ (y_{n}+1)^{2}}
\]
By (\ref{eq:xn}) we find that
\[
\lim_{n\rightarrow \infty }\frac{1-|\Psi (\alpha (\zeta _{n})
)|^{2}}{1-|\alpha (\zeta _{n})|^{2}}=\frac{1}{\ell_{\infty }} 
\]
Therefore, by Julia-Carath\'eodory's Theorem, $\Psi$ has a finite
angular derivative at $1$. 
More specifically, $\Psi^{\prime } (1)\leq
1/\ell_{\infty }$ (see Theorem \ref{thm:jc} (d)). Transfering back to
$\Hh $ we find that
\[
\frac{\psi (z)}{z}\longrightarrow A\geq \ell_{\infty }
\]
as $z$ tends to infinity non-tangentially. So, by Julia's Lemma,
$\psi (z)=Az+q (z)$,
where $\ima q (z)\geq 0$ on $\Hh$ and $q (z)/z$ tends to $0$ as $z$
tends to $\infty $ non-tangentially.
However,
\[
0\leq \ima q (\zeta _{n})=y_{n}-A\rightarrow\ell_{\infty }-A 
\]
So $A=\ell_{\infty }$.
\end{proof}

\begin{proof}[Proof of Corollary \ref{cor:conjtyp1}]
This follows the same lines argument as the proof of Theorem 2 of \cite{pom}.
By Theorem \ref{thm:conjtyp1}, we can write $\psi (z)=\ell_{\infty
}z+q (z)$, where $\ima q (z)\geq 0$ and $q (z)/z$ tends to $0$ as $z$
tends to $\infty $ non-tangentially. Set $U (z)=\ima q (z)/\ima
z$. Then by Schwarz-Pick on $q$, we find
\[
\left|\frac{q (w)-q (w^{\prime })}{w-w^{\prime }}\right|\leq \sqrt{U
(w)U (w^{\prime})}
\]
Let $L=\{t+i,t\geq 0 \}$. Then by (\ref{eq:latlimim}), $U (z)$ tends
to zero along $L$. On the other hand, from the Poisson integral
representaion of $\ima q$ it follows that $U$ is a decreasing function
of $y$. So given $\varepsilon >0$ there exist $R>0$ such that $U
(z)<\varepsilon $ for $z\in Q (R)$. Choose $\varepsilon $ much smaller
than $\ell_{\infty }$ and use 
\[
\left|\frac{\psi  (w)-\psi  (w^{\prime })}{w-w^{\prime }}\right|\geq
\ell_{\infty }-\left|\frac{q (w)-q (w^{\prime })}{w-w^{\prime
}}\right|  
\]
to conclude that $\psi $ is univalent on the corresponding quadrant $Q
(R)$.

The similar result about $\phi $ follows the exact same pattern,
therefore we omit the details.
\end{proof}

\begin{proof}[Proof of Proposition \ref{prop:philatlim}]
Recall from (\ref{eq:tauphitau}) that 
$H_{n} (z):=\tau _{n}^{-1}\circ \phi \circ
\tau _{n+1} (z)$ tends to $z$ uniformly on compact subsets of $\Hh $. 
Given a half-strip $\St (R,\ep )$, we can find a compact set $K$ such
that $\cup _{n=n_{0}}^{\infty }\tau _{n+1} (K)$ covers $\St (R^{\prime },\ep )$
for some $R^{\prime }>R$ and some $n_{0}$, because $y_{n}$ tends to
$\ell_{\infty }$ and $x_{n+1}-x_{n}$ tends to $\ell_{\infty }b_{0}$
(by Theorem \ref{thm:ufinheight}). A compution shows that,
\[
\tau _{n+1} (z)-\phi \circ \tau _{n+1} (z)=y_{n} (\tau _{n}^{-1}\circ \tau
_{n+1} (z)-H_{n} (z))\rightarrow  \ell_{\infty }b_{0}
\]
as $n$ tends to infinity, and this shows the proposition.
\end{proof}

\subsection{Non-existence results for parabolic maps of type I}\label{ssec:no}
\begin{proof}[Proof of Theorem \ref{thm:nonexist}]
Assume $\phi $  is a self-map of $\Hh $ of
parabolic type I.
Suppose that $w_{n}=x_{n}+iy_{n}$ is a BISBS. Then,
$z_{n}=\phi _{n} (w_{0})=u_{n}+iv_{n}$ has non-zero-step, and 
Pommerenke obtains a conjugation 
$\sigma $ with range in $\Hh$, see Theorem 
\ref{thm:forward}, where $\sigma \circ \phi =\sigma +b$ and $b\neq
0$ is the limit of $(u_{n+1}-u_{n})/v_{n}$, see (\ref{eq:b}). 
By Remark 1 of \cite{pom}, $\Arg z_{n}$ either tends to $0$ or to $\pi
$. Let us 
assume without loss of generality that it tends to $\pi$, so that
$b<0$.  Then, there exist $n_{0}$ so that for
$n\geq n_{0}$, $u_{n}-u_{n+1}$ is
positive, and 
greater than $(|b|/2) v_{n}$.
By telescoping 
sums
\[
|u_{n}|\geq (|b|/2)\sum _{k=n_{0}}^{n-1}v_{k} +|u_{n_{0}}|  \geq
(|b|/2) v_{n_{0}} 
(n-n_{0})+|u_{n_{0}}| 
\]
In particular,
\begin{equation}\label{eq:liminfun}
\liminf_{n\rightarrow \infty }\frac{|u_{n}|}{n}\geq (|b|/2) v_{n_{0}}
\end{equation}
This implies that we always have  $\liminf |u_{n}|/n>0$, 
and $\lim |u_{n}|/n=+\infty $ when $v_{n}\uparrow\infty $.

On the other hand, assume first that $\Arg w_{n}$ tends to $0$. Then,
$(x_{n+1}-x_{n})/y_{n}$ 
tends to $b_{0}>0$, see (\ref{eq:b0}).
So there exists $n_{1}$ such that for $n\geq n_{1}$, $x_{n+1}-x_{n}$
is less than $2b_{0}y_{n}$, and therefore
\[
x_{n}\leq 2b_{0}y_{n_{1}} (n-n_{1})+x_{n_{1}}
\]
A similar estimate can be obtained if $\Arg w_{n}$ tends to $\pi $. So
\begin{equation}\label{eq:limsupxn}
\limsup_{n\rightarrow \infty }\frac{|x_{n}|}{n}\leq 2b_{0}y_{n_{1}} 
\end{equation}
This implies that we always have $\limsup |x_{n}|/n<\infty $, 
and $\lim|x_{n}|/n=0 $ when $y_{n}\downarrow 0$.

Now, since $\sigma $ is a self-map of $\Hh $,
$d (\sigma  (w_{0}),\sigma (z_{n}))\leq d
(w_{0},z_{n})$ and the fact that $\sigma (z_{n})=\sigma (w_{0})+nb$ imply that
\[
\frac{v_{n}y_{0}}{(u_{n}-x_{0})^{2}+ (v_{n}+y_{0})^{2}}\leq 
\frac{(\ima \sigma (w_{0}))^{2}}{n^{2}b^{2}+4(\ima \sigma (w_{0}))^{2}}
\]
Now $v_{n}/u_{n}$ tends to zero, because $\Arg z_{n}$ tends to $\pi $. So
\begin{equation}\label{eq:SP}
\frac{v_{n}}{u_{n}^{2}}\leq \frac{C_{1}}{n^{2}}
\end{equation}
for some constant $C_{1}>0$.
On the other hand, let $\psi $ be the conjugation provided by Theorem
\ref{thm:conj}. The fact that $\psi (i-nb_{0})=z_{n}$ and
Schwarz-Pick's inequality $d (z_{n},w_{n})\leq d
(i-nb_{0},i+nb_{0})$ imply that
\[
\frac{1}{4n^{2}b_{0}^{2}+4}\leq 
\frac{v_{n}y_{n}}{(x_{n}-u_{n})^{2}+ (v_{n}+y_{n})^{2}}
\]
So that by (\ref{eq:SP}),
\begin{equation}\label{eq:bound}
y_{n} \geq  \frac{C_{2}}{n^{2}}\frac{(x_{n}-u_{n})^{2}}{v_{n}}
\geq   \frac{C_{2}}{C_{1}} \left(\frac{x_{n}}{u_{n}}-1\right)^{2}
\end{equation}
for some constant $C_{2}>0$. 

Assume now that $w_{n}$ is of type 2, i.e. $y_{n}\downarrow 0$.
Then (\ref{eq:bound}) yields
\begin{equation}\label{eq:xnun}
\frac{x_{n}}{u_{n}}\longrightarrow 1.
\end{equation}
But this contradicts the fact that in this case $\lim|x_{n}|/n=0 $ and $\liminf
|u_{n}|/n>0$. 
\end{proof}

\begin{proof}[Proof of Corollary \ref{cor:typeIa}]
Consider Pommerenke's map $\sigma $ (see
Theorem \ref{thm:forward}), in the case when $b\neq 0$ and the
imaginary parts of the forward iterates remain bounded above (type
Ia). Assume that $z_{n}=\phi _{n} (i)=u_{n}+iv_{n}$ and 
$\Arg z_{n}$ tends to $0$, so that $b>0$.
By Lemma 2 of \cite{pom}, $\sigma (u_{n}+zv_{n})-nb$ tends to $z$
uniformly on compact subsets of $\Hh $.  This can be used to show that
\begin{equation}\label{eq:latlimsig}
\latlim _{z\rightarrow +\infty }\rea \sigma (z)=+\infty 
\end{equation}
Suppose now that $\phi $ has a BISBS whose argument tends to
$0$ as well. By Theorem
\ref{thm:nonexist} such BISBS would have to be of type 1.
So it eventually belongs to some half-strip $\St (R,\ep)$.
But $\sigma (w_{n})=\sigma (w_{0})-nb$, so $\rea \sigma $ tends to
$-\infty $ along $w_{n}$, and this contradicts (\ref{eq:latlimsig}).
\end{proof}

\section{Examples}\label{sec:examples}

All of the examples that follow are obtained by conjugating $\phi $
to a translation on an appropriately chosen simply connected region,
so that the resulting Riemann map, which does the conjugation,  can be written
down explicitly.

\noindent {\bf Type Ia$\mathbf \emptyset$
(non-zero-step/finite-height/no-BISBS):}
\[
\phi (z)=\left[\sqrt{(\sqrt{z}+1)^{2}+1}-1 \right]^{2}
\]
Here $\phi $ is conjugated to translation by $1$ on the image of $\Hh
$ under the conformal map 
$(\sqrt{z}+1)^{2}$.
Then 
\[
z_{n}=\left(\sqrt{n+i}-1 \right)^{2}=n-1+i-2\sqrt{n+i}
\]
is a forward-iteration sequence, and there are no BISBS.

\noindent {\bf Type Ib$\mathbf \emptyset$
(non-zero-step/infinite-height/no-BISBS):}
\[
\phi (z)=\left(\sqrt{z}+1 \right)^{2}
\]
Here $\phi $ is conjugated to translation by $1$ on the image of $\Hh
$ under the conformal map $\sqrt{z}$.
Then 
\[
z_{n}= (n+i)^{2}=n^{2}-1+2in
\]
is a forward-iteration sequence, and there are no BISBS.
Note that $z_{n}/v_{n}= (n-1/n)+2i$ so it is non-zero-step.

\noindent {\bf Type IIb$\mathbf \emptyset$
(zero-step/infinite-height/no-BISBS):}
\[
\phi (z)=z+i
\]

\noindent {\bf Type Ia$\mathbf 1$
(non-zero-step/finite-height with a BISBS of non-zero-height):}

\[
\phi (z)=z+1
\]

\noindent {\bf Type IIb$\mathbf 1$
(zero-step/infinite-height with a BISBS of non-zero-height):}
\[
\phi (z)=\left(\sqrt{(\sqrt{z}-1)^{2}-1}+1 \right)^{2}
\]
Here $\phi $ is conjugated to translation by $-1$ on the image of $\Hh
$ under the conformal map $(\sqrt{z}-1)^{2}$.
Then
\[
z_{n}= (1+i\sqrt{n})^{2}=1-n+2i\sqrt{n}
\]
is a forward-iteration sequence, $z_{n}/v_{n}=-\sqrt{n}+1/\sqrt{n}+2i$
so the step goes to zero. Moreover,
\[
w_{n}=\left(\sqrt{n+i}+1 \right)^{2}=n+1+i+2\sqrt{n+i}
\]
is a BISBS of non-zero-height.

\noindent {\bf Type IIb$\mathbf 2$
(zero-step/infinite-height with a BISBS of zero-height):}
\[
\phi (z)=\sqrt{z^{2}-1}
\]
Here $\phi $ is conjugated to translation by $-1$ on the image of $\Hh
$ under the conformal map $z^{2}$.
Then 
\[
z_{n}=i\sqrt{n}
\]
is a forward orbit,
and
\[
w_{n}=\sqrt{n+i}
\]
is a BISBS of zero-height.

Examples of Ib1,IIa1,IIa2 type are still missing.

\end{document}